\documentclass[10pt,twoside]{article}
\usepackage{amsmath,amssymb,amsthm,amscd}
\usepackage{Latex-document}
\usepackage{latexsym}
\usepackage{float}
\usepackage{amssymb}
\usepackage{amsmath}
\usepackage[dvips]{graphicx}
\usepackage[dvips]{color}
\usepackage{colortbl}

\def\volumeno{I}
\title{\bf Statistical Equivalence and Stochastic\vskip -2mm Process Limit
Theorems\thanks{Research supported in part by NSF-Division of Mathematical Sciences.}\vskip 6mm}
\author{Lawrence D. Brown\thanks{Statistics Department,
Wharton School, University of Pennsylvania, Philadelphia, PA 19104-6340,
USA. E-mail: lbrown@wharton.upenn.edu}\vspace*{-0.5cm}}
\date{\vspace{-8mm}}

\begin{document}

\maketitle \thispagestyle{first} \setcounter{page}{557}

\begin{abstract}

\vskip 3mm

A classical limit theorem of stochastic process theory concerns the sample
cumulative distribution function (CDF) from independent random variables. If
the variables are uniformly distributed then these centered CDFs converge in
a suitable sense to the sample paths of a Brownian Bridge. The so-called
Hungarian construction of Komlos, Major and Tusnady provides a strong form
of this result. In this construction the CDFs and the Brownian Bridge sample
paths are coupled through an appropriate representation of each on the same
measurable space, and the convergence is uniform at a suitable rate.

Within the last decade several asymptotic statistical-equivalence theorems
for nonparametric problems have been proven, beginning with Brown and Low
(1996) and Nussbaum (1996). The approach here to statistical-equivalence is
firmly rooted within the asymptotic statistical theory created by L. Le Cam
but in some respects goes beyond earlier results.

This talk demonstrates the analogy between these results and those from the
coupling method for proving stochastic process limit theorems. These two
classes of theorems possess a strong inter-relationship, and technical
methods from each domain can profitably be employed in the other. Results in
a recent paper by Carter, Low, Zhang and myself will be described from this
perspective.

\end{abstract}

\vskip 12mm

\section{Probability setting}

\subsection{Background}

\vskip-5mm \hspace{5mm}

Let F be the CDF for a probability on [0,1];. F abs. cont., with
\[
f(x) \buildrel \Delta \over = {\frac{{\partial F}}{{\partial
x}}}\,\,on\,\,[0,1].
\]

Let X$_{1}$, \ldots , X$_{n}$ iid from F. $\hat {F}_{n} $ denotes the sample
CDF,
\[
\hat {F}_{n} (x) \buildrel \Delta \over = {\frac{{1}}{{n}}}{\sum\limits_{j =
1}^{n} {{\rm {\bf I}}_{[0,x]} \left( {X_{j}}  \right)}} .
\]

Let $\hat {Z}_{n} $denote the corresponding sample ``bridge'',
\begin{equation}
\hat {Z}_{n} (x) \buildrel \Delta \over = \hat {F}_{n} (x) - F(x)
\label{1}
\end{equation}

Let W(t) denote the standard Wiener process on [0,1] and let $\hat {W}_{n} $
denote the white noise process with drift f and local variance
${\raise0.7ex\hbox{${f(t)}$} \!\mathord{\left/ {\vphantom {{f(t)}
{n}}}\right.\kern-\nulldelimiterspace}\!\lower0.7ex\hbox{${n}$}}$. Thus
$\hat {W}_{n} $ solves
\[
d\hat {W}_{n} (t) = f(t)dt + \sqrt {{\frac{{f(t)}}{{n}}}} dW\left( {t}
\right).
\]

An alternate description of $\hat {W}_{n} $ is that it is the Gaussian
process with mean F(t) and independent increments having
\[
var\left( {\hat {W}_{n} (t) - \hat {W}_{n} (s)} \right) =
{\frac{{1}}{{n}}}\left( {F(t) - F(s)} \right),\,\,for\,\,0 \le s < t \le 1.
\]

The analog of $\hat {Z}_{n} $ is the Gaussian Bridge, defined by
\[
\hat {B}_{n} (t) = {\frac{{\hat {W}_{n} (t)}}{{\hat {W}_{n} (\ref{1})}}} -
F(t).
\]

There are various ways of describing the stochastic similarity between $\hat
{Z}_{n} $ and $\hat {B}_{n} $. For example Komlos, Major, and Tusnady (1975,
1976) proved a result of the form

\vspace{.1in}

\textbf{Theorem }(KMT): {\it Given any absolutely continuous F
{\{}X$_{1}$,\ldots ,X$_{n}${\}} can be defined on a probability
space on which $\hat {B}_{n} $ can also be defined as a
(randomized) function of {\{}X$_{1}$,\ldots ,X$_{n}${\}}. This can
be done in such a way that $\hat {B}_{n} $ has the Gaussian Bridge
distribution, above, and
\begin{equation}
P_{F} \left( {{\mathop {\sup} \limits_{t \in [0,1]}} \sqrt {n} {\left|
{\hat {Z}_{n} (t) - \hat {B}_{n} (t)} \right|}\, > a_{n}}  \right) \le c.
\label{2}
\end{equation}

Here $c>0$ and $a_n$ are suitable positive constants with $a_{n} \sim (d
\log n)/\sqrt n$ for some $d > 0$. The process $\hat {B}_{n} $ can be
constructed as a (randomized) function of $\hat {Z}_{n} $, that is,
$\hat {B}_{n} (t) = Q_{n} \left( {\hat {Z}_{n} (t)} \right)$. It should
be noted that the construction depends on knowledge of F.}

[Various authors, such as Cs\"{o}rg\"{o} and Revesz (1981) and Bretagnolle
and Massart (1989) have given increasingly detailed and precise values for
$a_n$ and $c=c(a_n)$, and also uniform (in n) versions of
(\ref{2}). These are not our focus.]

\subsection{Extensions}

\vskip-5mm \hspace{5mm}

1. Results like the above also extend to functional versions of the process
$\hat {Z}_{n} $. Various authors including Dudley (1978), Massart (1989) and
Koltchinskii (1994) have established results of the following form.

Let q:[0,1]$ \to \Re $ be of bounded variation. One can define
\[
\hat {Z}_{n} (q) \buildrel \Delta \over = \int {qd\{\hat {F}_{n}}  - F\} =
\int {\left( {F - \hat {F}_{n}}  \right)} dq.
\]

(Thus, $\hat {Z}_{n} (x) = \hat {Z}_{n} \left( {{\rm {\bf I}}_{[0,x]}}
\right)$.) There is a similar definition for $\hat {B}_{n} (q)$ as a
stochastic integral. (See, for example, Steele (2000).) Then the KMT theorem
extends to a fairly broad, but not universal, class of functions, $\cal Q$.
That is, for each F, $\hat {B}_{n} $ can be defined to satisfy
\begin{equation}
P_{F} \left( {{\mathop {\sup} \limits_{q \in {\rm Q}}} \sqrt {n}
{\left| {\hat {Z}_{n} (q) - \tilde {B}_{n} (q)\,} \right|} > {a}'_{n}}
\right) \le c \; {\mbox{ where}}\; {\{}{a}'_{n} {\}}\; \mbox{depends \;on}
\; \cal Q.
\label{3}
\end{equation}

(For most classes $\cal Q $, ${\raise0.7ex\hbox{${{a}'_{n} \sqrt {n}} $}
\!\mathord{\left/ {\vphantom {{{a}'_{n} \sqrt {n}}  {\log
n}}}\right.\kern-\nulldelimiterspace}\!\lower0.7ex\hbox{${\log n}$}} \to
\infty $ so that ${a}'_{n} > > a_{n} $.)
\vspace{.1in}

2. Bretagnolle and Massart (1989) proved a similar result for inhomogeneous
Poisson processes. Let {\{}T$_{1}$,\ldots ,T$_{N}${\}} be (ordered)
observations from an inhomogeneous Poisson process with cumulative intensity
function nF and, correpondingly, (local) intensity nf. Note that N$\sim
$Poisson(n) and conditionally given N the values of {\{}T$_{1}$,\ldots
,T$_{N}${\}} are the order statistics corresponding to an iid sample from
the distribution F. In this context we continue to define $\hat {F}_{n} (t)
= n^{ - 1}{\left\{ {{\sum\limits_{j = 1}^{N} {I_{[0,t]}}}  (T_{j} )}
\right\}}$ where the term in braces now has a Poisson distribution with mean
nF(t). Also, continue to define $\hat {Z}_{n} (t) \buildrel \Delta \over =
\hat {F}_{n} (t) - F(t)$ as in (\ref{1}). (But, note that it is no longer
true
that $\hat {Z}_{n} (\ref{1}) = 0$, w.p.1, as was the case in (\ref{1}).)

Then versions of the conclusions (\ref{2}) and (\ref{3}) remain valid. We
give an
explicit statement since this result will provide a model for our later
development.

\vspace{.1in} \textbf{Theorem }(BM): {\it Given any n and any
absolutely continuous F the observations {\{}T$_{1}$,\ldots
,T$_{N}${\}} of the inhomogeneous Poisson process can be defined
on a probability space on which $\hat {B}_{n} $ can also be
defined as a (randomized) function of {\{}T$_{1}$,\ldots
,T$_{N}${\}}. This can be done in such a way that $\hat {B}_{n} $
has the Gaussian Bridge distribution, above, and
\begin{equation}
P_{F} \left( {{\mathop {\sup} \limits_{t \in [0,1]}} \sqrt {n} {\left|
{\hat {Z}_{n} (t) - \hat {B}_{n} (t)} \right|} > a_{n}}  \right) \le c.
\label{4}
\end{equation}

Here $c>0$ and $a_n$ are suitable constants with $a_{n} \sim d \log n /
\sqrt n$.}

\textbf{Remark}: Clearly there must be extensions of (\ref{3}) that are
valid for
the Poisson case also, although we are not aware of an explicit treatment in
the literature. Such a statement would conclude in this setting that
\begin{equation}
P_{F} \left( {{\mathop {\sup} \limits_{q \in {\rm Q}}} \sqrt {n}
{\left| {\hat {Z}_{n} (q) - \hat {B}_{n} (q)} \right|}\, > {a}'_{n}}
\right) \le c \; {\mbox where} \;  {\{}{a}'_{n} {\}} \; \mbox{depends
\;on}\; \cal Q.
\label{5}
\end{equation}

\section{Main results}

\vskip-5mm \hspace{5mm}

The objective is a considerably modified version of (\ref{3}) and (\ref{5})
that is
stronger in several respects and (necessarily) different in others. We will
concentrate for most of the following on the statement (\ref{5}) since our
results
are slightly stronger and more natural in this setting. The extension of
(\ref{3})
will be deferred to a concluding Section.

Expression (\ref{5}) involves the target function $\hat {B}_{n} $. In the
modified
version the role of target function is instead played by $\tilde {W}_{n} $
which is the solution to the stochastic differential equation
\begin{equation}
d\tilde {W}_{n} (t) = g(t)dt + {\frac{{1}}{{2\sqrt {n}}} }dW(t)
\label{6}
\end{equation}
where $g(t) = \sqrt {f(t)} $. An alternate description of $\tilde
{W}_{n} $ is thus
\begin{equation}
\hspace{1in}  \tilde {W}_{n} = G(t) + W(t)/(2 \sqrt n)
\;{\mbox where}\; G(t) = {\int\limits_{0}^{t} {\sqrt {f(\tau )} d\tau}}.
\label{7}
\end{equation}
(In the special case where $f$ is the uniform density, $f$=1, then
$\tilde {W}_{n} = W_{4n} $.)

The role of the constructed random process $\hat {Z}_{n} $ is now played by
a differently constructed process $\tilde {Z}_{n} $. As before $\tilde
{Z}_{n} $ depends only on {\{}T$_{1}$,...T$_{N}${\}}, and not otherwise on
their CDF, F. This version also involves a large set, $\cal F $, of
absolutely
continuous CDFs. Both $\tilde {Z}_{n} $ and $\cal F $ will be described
later
in more detail. Here are statements of the main results.

\vspace{.1in} \textbf{Theorem 1}: {\it Let $\cal F $ be a set of
densities satisfying Assumption A or A', below. Let $\cal Q $ be the set
of all functions of bounded variation. Let {\{}T$_{1}$,\ldots
,T$_{N}${\}} be an inhomogeneous Poisson process with local intensity
$nf$. The process $\tilde {Z}_{n} $ can be constructed as a (randomized)
function of {\{}T$_{1}$,\ldots ,T$_{N}${\}}, with the construction not
depending on $f$. The Gaussian process $\tilde {W}_{n} $ having the
distribution (\ref{7} ) can also be defined on this same space as a
(randomized) function of {\{}T$_{1}$,\ldots ,T$_{N}${\}}. [This
construction depends on $f$ on a set of probability at most $c_{n}$.]
This can be done in such a way that
\begin{equation}
\sup _{f \in F} P_{f} \left( {{\mathop {\sup} \limits_{q \in {\rm Q}}
}{\left| {\tilde {Z}_{n} (q) - \tilde {W}_{n} (q)} \right|} > 0} \right) \le
c_{n} \rightarrow 0.
\label{8}
\end{equation}}

\vspace{.1in}

To be more precise, the phrase in brackets refers to the fact that there is
a basic construction, independent of $f$, and that this construction must
then
be modified on a set of measure at most $c_{n}$ with this set and the
modification depending on $f$.

For the situation of iid variables, as in (\ref{1}), a similar result holds.
In
this case the matching Gaussian process is again $\tilde {W}_{n} $, rather
than the Brownian bridge of the KMT theorem.

\vspace{.1in} \textbf{Theorem 2}: {\it Let $\cal F $ be a set of
densities satisfying Assumption B, below. Let $ \cal Q $ be the set of
all functions of bounded variation. Given any n and f$ \in \cal F $, iid
variables {\{}X$_{1}$,\ldots ,X$_{n}${\}} with density $f$ can be
defined on a probability space. A process $\tilde {\tilde {Z}}_{n} $ can
be constructed as a (randomized) function of {\{}X$_{1}$,\ldots
,X$_{n}${\}}, with the construction not depending on $f$. The Gaussian
process $\tilde {W}_{n} $ having the distribution (\ref{7}) can also be
defined on this same space as a (randomized) function of
{\{}X$_{1}$,\ldots ,X$_{n}${\}}. [This construction depends on $f$, but
only on a set of probability at most $c_{n}$.] This can be done in such
a way that
\begin{equation}
\sup _{f \in {\cal {F}}} P_{f} \left( {{\mathop {\sup} \limits_{q \in {\cal
{Q}}}
}{\left| {\tilde {\tilde {Z}}_{n} (q) - \tilde {W}_{n} (q)} \right|} > 0}
\right) \le c_{n} \to 0.
\label{9}
\end{equation}}

\section{Statistical background}

\subsection{Settings}

\vskip-5mm \hspace{5mm}

The first purpose of the discussion here is to motivate the probabilistic
results described above. A second purpose is to state the result on which to
base the proof of Theorem 1. The setting involves two statistical
formulations:

\underline {Formulation 1} (nonparametric inhomogeneous Poisson process):
The observations are \textbf{T} = {\{}T$_{1}$,\ldots ,T$_{N}${\}} from the
Poisson process with local intensity n$f$, $f \in \cal F $. The problem is
``nonparametric'' because the ``parameter space'', $\cal F $, is a very
large
set -- too large to be smoothly parameterized by a mapping from a (subset of)
a finite dimensional Euclidean space. Some possible forms for $\cal F $ are
discussed below. The statistician desires to make some sort of inference,
$\delta $, (possibly randomized) based on the observation of \textbf{X}.

\underline {Formulation 1'} (nonparametric density with random
sample size): The relation between Poisson processes and density problems
has been mentioned above. As a consequence, Problem 1 is equivalent to a
situation where the observations are {\{}X$_{1}$,...,X$_{N}${\}} with N$\sim
$Poisson(n) and {\{}X$_{1}$,...,X$_{N}${\}} the order statistics from a
sample of size N from the distribution with density $f$. Clearly, this
situation is closely related to the more familiar one in which the
observations are {\{}X$_{1}$,...,X$_{n}${\}} with n specified in advance.

\underline {Formulation 1''} (nonparametric density with
fixed sample size): This formulation refers to the more conventional density
setting in which the observations are {\{}X$_{1}$,...,X$_{n}${\}} iid with
density $f$.

\underline {Formulation 2} (white noise with drift): The statistician
observes a White noise process $d\tilde {W}_{n} (t)$, t$ \in $[0,1], with
drift g$ \in \cal G $ and local variance ${{\mbox{1}} \mathord{\left/
{\vphantom {{\mbox{1}} {\mbox{4n}}}} \right. \kern-\nulldelimiterspace}
{\mbox{4n}}}$. Thus
\[
d\tilde {W}_{n} (t) = g(t)dt + {\frac{{1}}{{2\sqrt {n}}} }dW(t),
\]
and $\tilde {W}_{n} (t) - G(t) = {\raise0.7ex\hbox{${W(t)}$}
\!\mathord{\left/ {\vphantom {{W(t)} {2\sqrt {n}
}}}\right.\kern-\nulldelimiterspace}\!\lower0.7ex\hbox{${2\sqrt {n}}
$}}$ where $G(t) = {\int\limits_{0}^{t} {g(\tau )d\tau}}  $. Again $\cal
G $ is a very large -- hence ``nonparametric'' -- parameter space.
Throughout, ${\cal G}   \subset {\cal{L}}_{2}= {\{}g: \smallint g^{2} <
\infty \}$. As of now, there need be no relation between $f$ in
Formulation 1 and g in Formulation 2, but such a relation will later be
assumed in connection with Theorem 1, where
\begin{equation}
g = \sqrt {f} \;\; {\mbox and} \; {\cal G} = {\left\{ {\sqrt {f} :f \in
{\cal F}}
\right\}}.
\label{10}
\end{equation}

This can alternatively be considered as a statistical formulation having
parameter space $\cal F $ under the identification (\ref{10}). We take this
point of
view in the BCLZ theorem, below.

\subsection{Constructive asymptotic statistical equivalence}

\vskip-5mm \hspace{5mm}

Here is one definition of the strongest form of such an equivalence.

\textbf{Definition} (asymptotic equivalence): {\it Let ${\cal P}
_{\mbox{j}}^{\mbox{(n)}} \mbox{ = (}{\cal
X}_{\mbox{j}}^{\mbox{(n)}} \mbox{, }{\cal
B}_{\mbox{j}}^{\mbox{(n)}} \mbox{,} {\cal
F}_{\mbox{j}}^{\mbox{(n)}} \mbox{)}$ j = 1,2, n = 1,2,... be two
sequences of statistical problems on the same sequence of
parameter spaces, $\Theta ^{(n)}$. Hence, ${\cal
F}_{\mbox{j}}^{\mbox{(n)}} = {\left\{ {F_{j,\theta} ^{(n)} :\theta
\in \Theta ^{(n)}} \right\}}$. Then $\Pi _{1}$ and $\Pi _{2} $are
asymptotically equivalent if there exist (randomized) mappings
$Q_{j}^{(n)} :{\cal X}_{\mbox{j}}^{\mbox{(n)}} \to {\cal
X}_{\mbox{k}}^{\mbox{(n)}} $, j, k = 1, 2, k$ \ne $j, such that
\begin{equation}
\sup _{\theta \in \Theta ^{(n)}} {\left\| {F_{j,\theta} ^{(n)} ( \cdot
) - \int {Q_{k}^{(n)} \left( {{\left. { \cdot}  \right|}x_{k}}
\right)F\left( {dx_{k}}  \right)}}  \right\|}_{TV} = c_{n} \to 0, j, k =
1, 2, k \ne j,
\label{11}
\end{equation}
where ${\left\| { \cdot}  \right\|}_{TV} $ denotes the total
variation norm.}

This definition involves a reformulation of the general theory originated by
LeCam (1953, 1964). See also Le Cam (1986), Le Cam and Yang (2000), van der
Vaart (2002) and Brown and Low (1996) for background on this theory
including several alternate versions of the definition and related concepts,
a number of conditions that imply asymptotic equivalence, and many
applications to a variety of statistical settings. Note that both
Formulations 1 and 2 involve an index, n, and can thus be considered as
sequences of statistical problems in the sense of the definition.

\subsection{Spaces of densities (or intensities)}

\vskip-5mm \hspace{5mm}

Suitable families of densities, $\cal F $, can be defined via Besov norms
with
respect to the Haar basis. The Besov norm with index $\alpha $ and shape
parameters p = q can most conveniently be defined via the stepwise
approximants to $f$ at resolution level k. These approximants are defined as
\[
\bar {f}_{k} (t) = {\sum\limits_{\ell = 0}^{2^{k} - 1} {I_{[{{\ell}
\mathord{\left/ {\vphantom {{\ell}  {2^{k}}}} \right.
\kern-\nulldelimiterspace} {2^{k}}},{{(\ell + 1)} \mathord{\left/ {\vphantom
{{(\ell + 1)} {2^{k}}}} \right. \kern-\nulldelimiterspace} {2^{k}}})}
(t){\int\limits_{{{\ell}  \mathord{\left/ {\vphantom {{\ell}  {2^{k}}}}
\right. \kern-\nulldelimiterspace} {2^{k}}}}^{{{(\ell + 1)} \mathord{\left/
{\vphantom {{(\ell + 1)} {2^{k}}}} \right. \kern-\nulldelimiterspace}
{2^{k}}}} {2^{k}f}}} } ,
\]
and the Besov($\alpha $,p) norm is defined as
\[
{\left\| {f} \right\|}_{\alpha ,p} = {\left\{ {{\left| {\bar {f}_{0}}
\right|}^{p} + {\sum\limits_{k = 0}^{\infty}  {2^{pk\alpha} {\left\| {\bar
{f}_{k} - \bar {f}_{k + 1}}  \right\|}_{p}^{p}}} }  \right\}}^{{{1}
\mathord{\left/ {\vphantom {{1} {p}}} \right. \kern-\nulldelimiterspace}
{p}}}.
\]

The statement of Theorem 1 can now be completed by stating the assumption on
$\cal F $ needed for its validity.

\textbf{Assumption A}: {\it $\cal F $ satisfies
\begin{equation}
{\cal F} \subset {\left\{ {f:\inf _{0 \le x \le 1} f(x) \ge \varepsilon _{0}
} \right\}} \; {\mbox  for \; some}\;  \varepsilon _{0}>0
\label{12}
\end{equation}
and $\cal F $ is compact in both Besov(1/2,2) and Besov(1/2,4).}

Other function spaces are also conventional for nonparametric statistical
applications of this type. The most common of these are based on either the
Lipshitz norm ${\left\| {f} \right\|}_{\beta} ^{(L)} $ or the Sobolev norm
${\left\| {f} \right\|}_{\beta} ^{(S)} $. These are defined for $\beta
\le $1 by
\[
{\left\| {f} \right\|}_{\beta} ^{(L)} = {\mathop {\sup} \limits_{0 \le x < y
\le 1}} {\frac{{{\left| {f(y) - f(x)} \right|}}}{{{\left| {y - x}
\right|}^{\beta}} }},\,\,\,\,\,\,{\left\| {f} \right\|}_{\beta} ^{(S)} =
{\sum\limits_{ - \infty} ^{\infty}  {k^{2\beta} \vartheta _{k}^{2}}}
\]
where $\vartheta _{k} = {\int\limits_{0}^{1} {f(x)e^{ik2\pi x}dx}} $
denote the Fourier coefficients of f. (Both spaces have natural
definitions for $\beta $>1 as well, but we need consider here only the
case $\beta  \le 1$.)

The following implies Assumption A and hence also suffices for validity of
Theorem 1.

\textbf{Assumption A':} {\it $\cal F$ satisfies (\ref{12}), and is
bounded in the Lipshitz norm with index $\beta $, and is compact
in the Sobolev norm with index $\alpha $, where $\alpha   \ge
\beta $ and either $\beta  > 1/2$ or $\alpha  \ge  3/4$ and
$\alpha +\beta  \ge 1$.}

The following assumption is noticeably stronger than either A' or A, and is
used in Theorem 2.

\textbf{Assumption B}: {\it $\cal F $ satisfies (\ref{12}) and is
bounded in the Lipshitz norm with index $\beta $, where $\beta >
1/2$.}

For more information about the relation of these spaces in this context see
Brown, Cai, Low and Zhang (2002) and Brown, Carter, Low and Zhang (2002)
(referred to as BCLZ below).

\subsection{Statistical equivalence theorems}

\vskip-5mm \hspace{5mm}

BCLZ then extended earlier results of Nussbaum (1996) and Klemela and
Nussbaum (1998) to prove the following basic result:

\vspace{.1in} \textbf{Theorem a} (BCLZ): {\it Consider the statistical
Formulations 1 and 2 with the parameter space ${\cal F} $ and the
relation (\ref{10}). Assume ${\cal F} $ satisfies Assumption A (or A'.
Then the sequences of statistical problems defined in these two
formulations are asymptotically statistically equivalent.}

BCLZ describes in detail a construction of $\tilde {Z}_{n} $ as a
(randomized) function of {\{}T$_{1}$,\ldots ,T$_{n}${\}}. (More precisely,
BCLZ describes the construction of the Haar basis representation of $\tilde
{Z}_{n} $, from which $\tilde {Z}_{n} $ can directly be recovered.) This
construction is invertible, in that {\{}T$_{1}$,\ldots ,T$_{n}${\}} can be
recovered as a function of $\tilde {Z}_{n} $. Further, BCLZ shows that both
$\tilde {Z}_{n} $and $\tilde {W}_{n} $ can be represented on the same
probability space so that their distributions, $P_{\tilde {Z}_{n}}  $ and
$P_{\tilde {W}_{n}}  $, say, satisfy
\[
{\left\| {P_{\tilde {Z}_{n}}  - P_{\tilde {W}_{n}}}   \right\|}_{TV} \to 0.
\]

The mappings {\{}$Q_{j}^{(n)} $: j=1,2, n = 1,2,...{\}} that yield the
equivalence of the above theorem can then be directly inferred from this
construction. To save space here we refer the reader to that paper or Brown (2002) for
details of the construction and proof. It can be remarked that these bear
considerable similarity to parts of the construction and proof in
Bretagnolle and Massart (1989) and other proofs of KMT type theorems. But
there are also some basic differences, especially those related to the
appearance of the square-root in the fundamental relation (\ref{10}) and the
total
variation norm in the definition of equivalence. In addition, the fact that
(\ref{8}) is uniform in $ \cal Q $ and ${\cal F} $ entails the need for
various refinements
in the proof.

Theorem 1 is now an immediate logical consequence of this result from BCLZ
and the following lemma.

\vspace{.1in} \textbf{Lemma}: {\it Suppose ${\cal P}
_{\mbox{j}}^{\mbox{(n)}} \mbox{ = (}{\cal
X}_{\mbox{j}}^{\mbox{(n)}} \mbox{,} {\cal
B}_{\mbox{j}}^{\mbox{(n)}} \mbox{, }{\cal
F}_{\mbox{j}}^{\mbox{(n)}} \mbox{)}$ j = 1,2, n = 1,2,... are
asymptotically equivalent sequences of statistical problems on the
same sequence of parameter spaces, $\Theta ^{(n)}$. Let
{\{}$Q_{j}^{(n)} $: j=1,2, n=1,2,...{\}} denote a sequence of
mappings that define this equivalence, as in (\ref{11}). Then
there are non-randomized mappings {\{}$\tilde {Q}_{j}^{(n)} $:
j=1,2, n = 1,2,...{\}} such that
\begin{equation}
P_{f} (\tilde {Q}_{j}^{(n)} = Q_{j}^{(n)} ) \ge 1 - c_{n} \; {\mbox for \;
every }\;
f \in {\cal F}_{\mbox{j}}^{\mbox{(n)}}, \; j=1, 2, n=1, 2,...
\label{13}
\end{equation}
and for every $\theta  \in \Theta ^{(n)}$
\begin{equation}
P_{f_{j,\theta}}   \left( {\tilde {Q}_{j}^{(n)} \left( {X_{j}^{(n)}}
\right) \in A} \right) = P_{f_{k,\theta}}   \left( {X_{k}^{(n)} \in A}
\right), \theta  \in \Theta ^{(n)}
\label{14}
\end{equation}
for every measurable $A \subset {\cal X}_{\mbox{k}}^{\mbox{(n)}} $, $j
,k=1, 2$, $j \ne k$, $n=1,2,\cdots$.}

\textbf{Proof of Lemma}: Fix n, j, k$ \ne $j, $\theta  \in \Theta
^{(n)}$. Let $F_{k}$ denote the distribution under $\theta $ of
$X_{k}^{(n)} $ and let ${F}'_{k} $ denote the distribution under $\theta $
of $Q_{j}^{(n)} \left( {X_{j}^{(n)}}  \right)$. Let H = min($F_{k}$,
${F}'_{k} )$. Let $\infty \ge {f}'_{k} = {\frac{{d{F}'_{k}}} {{dH}}} \ge 1$.
Then define $\tilde {Q}_{j}^{(n)} $ as a version of the randomized map
satisfying
\[
\tilde {Q}_{j} (B\vert x) = {\frac{{1}}{{{f}'_{k} (x)}}}Q_{j} (B\vert x)\, +
\,{\frac{{{f}'_{k} - 1}}{{{f}'_{k}}} }\left( {{F}'_{k} (B) - H(B)} \right).
\]

This completes the proof of the lemma, and consequently also that of Theorem
1. $\Box$

\vspace{.1in}
Theorem 2 requires a slightly different fundamental result. The following
result is the foundation for the proof of Theorem 2. It is adapted from
Theorem 2 of BCLZ. This result closely resembles Theorem a, above, but as
noted in BCLZ it appears to require a modified construction for its proof.
The argument there is based heavily on results in Carter (2001).

\vspace{.1in} \textbf{Theorem b} (BCLZ): {\it Consider the
statistical Formulations 1'' and 2 with the parameter space ${\cal
F} $ and the relation (\ref{10}). Assume ${\cal F} $ satisfies
Assumption B. Then the sequences of statistical problems defined
in these two formulations are asymptotically statistically
equivalent.}

\label{lastpage}

\end{document}